\documentclass[12pt]{amsart}
\usepackage{amscd,amsmath,amssymb,amsfonts}
\usepackage[cmtip, all]{xy}
\begin{document}
\newtheorem{thm}[subsection]{Theorem}
\newtheorem{prop}[subsection]{Proposition}
\newtheorem{lem}[subsection]{Lemma}
\newtheorem{cor}[subsection]{Corollary}
\theoremstyle{definition}
\newtheorem{Def}[subsection]{Definition}
\theoremstyle{remark}
\newtheorem{rem}[subsection]{Remark}
\newtheorem{rems}[subsection]{Remarks}
\newtheorem{ex}[subsection]{Example}
\newtheorem{exs}[subsection]{Examples}
\numberwithin{equation}{section}
\newcommand{\an}{{\rm an}}
\newcommand{\alg}{{\rm alg}}
\newcommand{\cl}{{\rm cl}}
\newcommand{\Alb}{{\rm Alb}}
\newcommand{\CH}{{\rm CH}}
\newcommand{\mc}{\mathcal}
\newcommand{\mb}{\mathbb}
\newcommand{\surj}{\twoheadrightarrow}
\newcommand{\inj}{\hookrightarrow}
\newcommand{\red}{{\rm red}}
\newcommand{\codim}{{\rm codim}}
\newcommand{\rank}{{\rm rank}}
\newcommand{\Pic}{{\rm Pic}}
\newcommand{\Div}{{\rm Div}}
\newcommand{\Hom}{{\rm Hom}}
\newcommand{\im}{{\rm im}}
\newcommand{\Spec}{{\rm Spec \,}}
\newcommand{\sing}{{\rm sing}}
\newcommand{\Char}{{\rm char}}
\newcommand{\Tr}{{\rm Tr}}
\newcommand{\Gal}{{\rm Gal}}
\newcommand{\Min}{{\rm Min \ }}
\newcommand{\Max}{{\rm Max \ }}
\newcommand{\supp}{{\rm supp}\,}
\newcommand{\0}{\emptyset}
\newcommand{\sA}{{\mathcal A}}
\newcommand{\sB}{{\mathcal B}}
\newcommand{\sC}{{\mathcal C}}
\newcommand{\sD}{{\mathcal D}}
\newcommand{\sE}{{\mathcal E}}
\newcommand{\sF}{{\mathcal F}}
\newcommand{\sG}{{\mathcal G}}
\newcommand{\sH}{{\mathcal H}}
\newcommand{\sI}{{\mathcal I}}
\newcommand{\sJ}{{\mathcal J}}
\newcommand{\sK}{{\mathcal K}}
\newcommand{\sL}{{\mathcal L}}
\newcommand{\sM}{{\mathcal M}}
\newcommand{\sN}{{\mathcal N}}
\newcommand{\sO}{{\mathcal O}}
\newcommand{\sP}{{\mathcal P}}
\newcommand{\sQ}{{\mathcal Q}}
\newcommand{\sR}{{\mathcal R}}
\newcommand{\sS}{{\mathcal S}}
\newcommand{\sT}{{\mathcal T}}
\newcommand{\sU}{{\mathcal U}}
\newcommand{\sV}{{\mathcal V}}
\newcommand{\sW}{{\mathcal W}}
\newcommand{\sX}{{\mathcal X}}
\newcommand{\sY}{{\mathcal Y}}
\newcommand{\sZ}{{\mathcal Z}}
\newcommand{\A}{{\mathbb A}}
\newcommand{\B}{{\mathbb B}}
\newcommand{\C}{{\mathbb C}}
\newcommand{\D}{{\mathbb D}}
\newcommand{\E}{{\mathbb E}}
\newcommand{\F}{{\mathbb F}}
\newcommand{\G}{{\mathbb G}}
\renewcommand{\H}{{\mathbb H}}
\newcommand{\I}{{\mathbb I}}
\newcommand{\J}{{\mathbb J}}
\newcommand{\M}{{\mathbb M}}
\newcommand{\N}{{\mathbb N}}
\renewcommand{\P}{{\mathbb P}}
\newcommand{\Q}{{\mathbb Q}}
\newcommand{\R}{{\mathbb R}}
\newcommand{\T}{{\mathbb T}}
\newcommand{\U}{{\mathbb U}}
\newcommand{\V}{{\mathbb V}}
\newcommand{\W}{{\mathbb W}}
\newcommand{\X}{{\mathbb X}}
\newcommand{\Y}{{\mathbb Y}}
\newcommand{\Z}{{\mathbb Z}}

\newcommand{\Nm}{{\operatorname{Nm}}}
\newcommand{\NS}{{\operatorname{NS}}}
\newcommand{\id}{{\operatorname{id}}}

\title{The Steinberg Curve}

\author{H\'el\`ene Esnault
\and Marc Levine}

\address{
Universit\"{a}t GH Essen\\
FB6 Mathematik und Informatik\\
45117 Essen\\
Germany\\
\and
Department of Mathematics\\
Northeastern University\\
Boston, MA 02115\\
USA}

\thanks{Partially supported by the NSF and by the
DFG-Forschergruppe ``Arithmetik und Geometrie''}

\email{ esnault@uni-essen.de marc@neu.edu}



\renewcommand{\abstractname}{Abstract}
\begin{abstract}
Let $E$ and $E'$ be elliptic curves over $\C$, with Tate parametrizations 
$p:\C^*\to E$, $p':\C^*\to E'$. We have the map $p*p':\C^*\otimes\C^*\to 
F^2 \CH^2(E\times E')$
sending $u\otimes v$ to the class of the zero cycle 
$(x,y)-(x,0)-(0,y)+(0,0)$, where
$x=p(u)$, $y=p'(v)$. We show that, for general $u\in\C^*$, 
$p*p'(u\otimes(1-u))$ is
not zero in $\CH^2$. We also show that  the cycle $p*p'(u\otimes(1-u))$ 
is not
detectable by a certain class of cohomology theories, including the 
cohomology of the analytic motivic complex involving the dilogarithm 
function
defined by S. Bloch in \cite{Bl}. This is in contrast to its
\'etale version defined by S. Lichtenbaum \cite{Li}, which contains
the Chow group. \ \\
\end{abstract}

\maketitle

\section{Tate curves and line bundles}
For a scheme $X$ over $\C$, we let $X_\an$ denote the set of $\C$-points 
with the
classical topology. We let $\sO_{X_\an}$ denote the sheaf of holomorphic 
functions on
$X_\an$.

We begin by describing a construction of the universal analytic Tate 
curve over 
$\C$. We first form
the analytic manifold $\hat\sC^*$ as 
the quotient of the disjoint union $\sqcup_{i=-\infty}^\infty U_i$, 
with each $U_i=\C^2$, by the equivalence relation
$$(x,y)\in U_i\setminus\{Y=0\} \sim (\frac{1}{y},xy^2)\in 
U_{i+1}\setminus\{X=0\}.$$
The function $\tilde\pi(x,y)=xy$ on $\hat\sC^*$ is globally
defined. Letting 
$D\subset\C$
be the disk $\{|z|<1\}$, we define $\sC^*=\tilde\pi^{-1}(D)$, so 
$\tilde\pi$ 
restricts
to the analytic map
$\pi:\sC^*\to D$. We let $\tilde 0:D\to\sC^*$ be the section $z\mapsto 
(z,1)\in U_0$.

Let $D^*\subset D$ be the punctured disk $z\ne0$. Since the map
$(x,y)\mapsto(\frac{1}{y},xy^2)$ is an automorphism of
$(\C^*)^2$, the open submanifold
$\pi^{-1}(D^*)$ of $\sC^*$ is isomorphic to $(\C^*)^2$, and the
restriction of the map $\pi$ is just the map 
$(x,y)\mapsto xy$. 
Thus, the projection $p_2:(\C^*)^2\to\C^*$ gives an
isomorphism of the fiber $\sC^*_t:=\pi^{-1}(t)$ with
$\C^*$, for $t\in D^*$. 

The fiber $\pi^{-1}(0)$, on the other hand, is an
infinite union of projective lines. Indeed, define the map $f_i:\C\P^1\to 
\sC^*_0$
by sending $(a:1)\in \C\P^1\setminus\infty$ to
$(0,a)\in U_i$, and $\infty=(1:0)$ to $(0,0)\in U_{i+1}$, and let 
$C_i=f_i(\C\P^1)$.
Then $\pi^{-1}(0)=\cup_{i=-\infty}^\infty C_i$, with $\infty\in C_i$ 
joined
with $0\in C_{i+1}$. Note in particular that the value
$\tilde{0}(0)$ of the zero section avoids the singularities of 
$\pi^{-1}(0)$.

Define the automorphism $\phi$ of $\sC^*$ over $D$ by sending $(x,y)\in 
U_i$ to
$(x,y)\in U_{i-1}$. This gives the action of $\Z$ on $\sC^*$, with $n$ 
acting by
$\phi^n$. It is easy to see that this action is free and proper, so the 
quotient space
$\sE:=\sC^*/\Z$ exists as a bundle $\pi:\sE\to D$. The section $\tilde 
0:D\to\sC^*$
induces the section $0: D\to\sE$.

Take $t\in D^*$. Identifying $\sC^*_t$ with $\C^*$ as above, we see that 
$\phi$
restricts to the automorphism $z\mapsto tz$. Thus, the fiber 
$\sE_t:=\pi^{-1}(t)$ for
$t\in D^*$ is the {\em Tate elliptic curve} $\C^*/t^\Z$, with identity 
$0(t)$. On
$\sC^*_0$, however, $\phi$ is the union of the ``identity" isomorphisms 
$C_i\to
C_{i-1}$. Thus $\phi(\infty\in C_i)=0\in C_i$, so the restriction of 
$\sC^*_0\to
\sE_0$ to $C_0$ identifies $\sE_0$  with the nodal curve $\C\P^1/0\sim 
\infty$.  We
let $*\in \sE_0$ denote the singular point. Then $\tilde{0}(0)
\in \sE_0 \setminus *$.

The map $(t,w)\in D\times\C^*\mapsto (\frac{t}{w},w)\in U_0$ gives an 
isomorphism  
$\psi:D\times\C^*\to U_0\setminus\{Y=0\}$ over $D$. The composition
\[
D\times\C^*\to U_0\setminus\{Y=0\}\subset\sC^*\xrightarrow{q}\sE
\]
defines the map $p:D\times\C^*\to \sE$ over $D$, with image 
$\sE\setminus\{*\}$. 

Take $u\in\C^*$. We have the local system on $\sE$
\[
\sL_u:= \sC^*\times\C/(z,\lambda)\sim(\phi(z),u\lambda)\to\sE,
\]
and the associated holomorphic line bundle $\sL^\an_u$ on $\sE$. 

Let $E_t$ be the algebraic elliptic curve associated to the analytic 
variety $\sE_t$,
let $L_u(t)$ and $L_u^\an(t)$ denote the restriction of $\sL_u$ and 
$\sL_u^\an$ to
$\sE_t$, and let $L_u^\alg(t)$ be the algebraic line bundle on $E_t$ 
corresponding to
$L_u^\an(t)$ via \cite{Se}. The restriction of $p$ to $t\times\C^*$ 
defines the map
$p_t:\C^*\to E_{t\an}$. For $t\ne0$, $p_t$ is a covering space of 
$E_{t\an}$. The map 
$p_0:\C^*\to E_{0\an}$ is the analytic map associated to the algebraic 
open immersion
\[
\P^1\setminus\{0,\infty\}\xrightarrow{j}\P^1\to \P^1/0\sim\infty=E_0.
\]

If $E$ be an elliptic curve over $\C$, then $E_\an\cong\C/\Lambda$, where
$\Lambda\subset\C$ is a lattice spanned by $1$ and some $\tau$ in the 
upper half
plane. Taking $t=e^{2\pi i\tau}$ gives the isomorphism $E_\an\cong\sE_t$, 
so each
elliptic curve over $\C$ occurs as an $E_t$ for some (in fact for 
infinitely many)
$t\in D^*$.

Sending $u\in\C^*$ to the isomorphism class of $L_u^\alg(t)$ defines 
a homomorphism $\tilde p_t:\C^*\to\Pic(E_t)$. We denote the identity 
$0(t)\in E_t$
simply by 0 if $t$ is given.

\begin{lem}\label{C1Comp} For all $t\in D$,  
$c_1(L_u^\alg(t))=(p_t(u))-(0)$.
\end{lem}

\begin{proof} We first handle the case $t\ne0$. Let $q:\C\to E:=E_t$ be 
the map
$q(z)=p_t(e^{2\pi iz})$, let
$\tau\in
\C$ be an element with $e^{2\pi i\tau}=t$, and let $\Lambda\subset\C$ be 
the lattice
generated by $1$ and $\tau$. The map $q$ identifies
$E$ with $\C/\Lambda$, and $L_u(t)$ with the local system defined by the 
homomorphism
$\rho:\Lambda\to\C^*$, $\rho(a+b\tau)=u(b)$.

There is a unique cocycle $\theta$ in 
$Z^1(\Lambda,H^0(\C,\sO^*_{\C_\an}))$ with
$\theta(1)=1$, $\theta(\tau)=e^{-2\pi iz}$; let $L$ be the corresponding 
holomorphic
line bundle on $E$. Computing $c_1^{\rm top}(L)\in H^2(E,\Z)$ by using 
the 
exponential
sequence, we find that $\deg(L)=1$. By Riemann-Roch, we have 
$H^0(E,L)=\C$; let
$\Theta(z)$ be the corresponding global holomorphic function on $\C$, 
i.e.,
\[
\Theta(z+1)=\Theta(z),\ \Theta(z+\tau)=e^{-2\pi iz}\Theta(z),
\]
and the divisor of $\Theta$ on $E$ is $(x)$, with $L\cong\sO_E(x)$.

Take $v,w\in\C$ with $u=e^{2\pi iv}$ and $q(w)=x$.  Let
$f(z)=\frac{\Theta(z+w-v)}{\Theta(z+w)}$.
Then
\[
f(z+1)=f(z),\  f(z+\tau)=uf(z),
\]
and $\Div(f)=(p(u))-(0)$. Thus, multiplication by $f$ defines an
isomorphism
\[
\times f:\sO_{E_\an}((p(u))-(0))\to L_u^\an.
\]

The proof for $E_0=\P^1/0\sim\infty$ is essentially the same, where we 
replace 
$\frac{\Theta(z+w-v)}{\Theta(z+w)}$ with the rational function
$\frac{X-u}{X-1}$.
\end{proof}

Thus, the image of $\tilde p_t$ in $\Pic(E_t)$ is $\Pic^0(E)$. After 
identifying the
smooth locus of $E_t^0$ of
$E_t$ with 
$\Pic^0(E_t)$
by sending $x\in E_t^0$ to the class of the invertible sheaf 
$\sO_{E_t}((x)-(0))$, we have
$\tilde p_t=p_t$.

\section{The Albanese kernel and the Steinberg relation}\label{AlbKer}
Let $X$ be a smooth projective variety. We let $\CH_0(X)$ denote the 
group of zero
cycles on $X$, modulo rational equivalence, $F^1\CH_0(X)$ the subgroup of 
cycles of
degree zero, and $F^2\CH_0(X)$ the kernel of the Albanese map
$\alpha_X:F^1\CH_0(X)\to\Alb(X)$. The choice of a point $0\in X$ gives a 
splitting to
the inclusion $F^1\CH_0(X)\to \CH_0(X)$. 

Let $E$, $E'$ be smooth elliptic curves.  As $\Alb(E\times E')=E\times 
E'$, the 
inclusion
$F^2\CH_0(E\times E')\to F^1\CH_0(E\times E')$ is split by sending 
$(x,y)-(0,0)$ to 
$(x,y)-(x,0)-(0,y)+(0,0)$. Thus $F^2\CH_0(E\times E')$ is generated by 
zero-cycles
of the form $(x,y)-(x,0)-(0,y)+(0,0)$. Choosing an isomorphism $E\cong 
E_t$, $E'\cong
E_{t'}$, we have the covering spaces $p:\C^*\to E_\an$, $p':\C^*\to 
E'_\an$, and the  map
\begin{gather}
p*p':\C^*\otimes\C^* \to F^2\CH_0(E\times E')\label{BasicSurjDef1}\\
u\otimes v \mapsto 
p(u)*p'(v):=\notag \\
(p(u),p'(v))-(p(u),0)-(0,p'(v))+(0,0).\notag
\end{gather}
By the theorem of the cube \cite{M}, the map $p*p'$ is a group
homomorphism, and thus is surjective.

In case one or both of $E$, $E'$ is the singular curve $E_0$, we will 
need to use the
the theory of zero-cycles mod rational equivalence defined in 
\cite{LevWeib}.  If $X$
is a reduced, quasi-projective variety over a field $k$ with singular 
locus
$X_\sing$, the group $\CH_0(X)$ (denoted $\CH_0(X, X_\sing)$ in 
\cite{LevWeib})  is
defined as the quotient of the free abelian group on the regular closed 
points of
$X$, modulo the subgroup generated by zero-cycles of the form $\Div f$, 
where $f$ is
a rational function on a dimension one closed subscheme
$D$ of $X$ such that
\begin{enumerate}
\item No irreducible component of $D$ is contained in $X_\sing$.
\item In a neighborhood of each point of $D\cap X_\sing$, the subscheme
$D$ is a complete intersection.
\item $f$ is in the subgroup $\sO_{D,D\cap X_\sing}^*$ of $k(D)^*$.
\end{enumerate}
It follows in particular 
from these conditions that $\Div f$ is a sum of regular points 
of $X$.

For $X$ a reduced curve, sending a regular closed point $x\in X$ to the 
invertible
sheaf $\sO_X(x)$ extends to give an isomorphism $\CH_0(X)\cong \Pic(X)$.

We  extend the definition of $F^i\CH_0$ to $E\times E'$ with either 
$E=E_0$ or
$E'=E_0$ or
$E=E'=E_0$, by defining  $F^1\CH_0(E\times E')$ as the 
subgroup of $\CH_0(E\times E')$ generated by the differences $[x]-[y]$, 
and
$F^2\CH_0(E\times E')$ the  subgroup generated by expressions
$[(x,y)]-[(x,0)]-[(0,y)]+[(0,0)]$, where $x$ is a smooth point
of $E$ and $y$ a smooth point of $E'$. The surjection  
$p*p':\C^*\otimes\C^*\to
F^2\CH_0(E\times E')$ is then defined by the same formula as 
\eqref{BasicSurjDef1}.

\begin{prop}[The Steinberg relation] Take $E=E'=E_0$. Then 
$p(u)*p(1-u)=0$ in
$\CH_0(E_0\times E_0)$ for all $u\in \C\setminus\{0,1\}$.
\end{prop}

\begin{proof} Let $X$ be a quasi-projective surface over a field $k$. By 
\cite{Lev}, there is an isomorphism $\phi:H^2(X,\sK_2)\to CH_0(X)$. The 
product
$\sO_X^*\otimes\sO_X^*\to\sK_2$ gives the cup product 
\[
H^1(X,\sO_X^*)\otimes H^1(X,\sO_X^*)\xrightarrow{\cup} H^2(X,\sK_2).
\]
In addition, let $D$, $D'$ be Cartier divisors which intersect 
properly on $X$, and suppose that $\supp D\cap\supp D'\cap X_\sing=\0$. 
Then
\begin{equation}\label{ProdComp}
\phi(\sO_X(D)\cup\sO_X(D'))=[D\cdot D'],
\end{equation}
where $\cdot$ is the intersection product and $[-]$ denotes the class in 
$\CH_0$.

Since  $L_u^\alg=\sO_{E_0}(p(u)-0)$, \eqref{ProdComp} implies
\[
p(u)*p(1-u)=\rho(p_1^*L_u^\alg\cup p_2^*L_{1-u}^\alg),
\]
so it suffices to show that $p_1^*L_u^\alg\cup p_2^*L_{1-u}^\alg=0$ in
$H^2(E_0\times E_0,\sK_2)$.

Write $X$ for $E_0\times E_0$. Let $\bar \sK_2$ be the image of $\sK_2$ 
in the
constant sheaf $K_2(\C(X))$. By Gersten's conjecture, the surjection
$\pi:\sK_2\to\bar\sK_2$ is an isomorphism at each regular point of $X$, 
hence $\pi$
induces an isomorphism on $H^2$.

Let $q:\P^1\to E_0$ be the normalization, giving the normalization 
$q\times
q:\P^1\times\P^1\to X$. Let $i:*\to E_0$ be the inclusion of the singular 
point. We
have the exact sequence of sheaves on $E_0$
\begin{equation}\label{K1Comp}
q_*\sK_1\xrightarrow{\beta} i_*K_1(\C)\to0
\end{equation}
and the exact sequence of sheaves on $X$:
\begin{equation}\label{K2Comp}
(q\times q)_*\sK_2\xrightarrow{\alpha} (i\times
q)_*\sK_2\oplus(q\times i)_*\sK_2
\to (i\times i)_* K_2(\C)\to 0,
\end{equation} 
with augmentations $\epsilon_1:\sK_1\to\eqref{K1Comp}$, 
$\epsilon_2:\bar\sK_1\to\eqref{K2Comp}$. The various cup products in 
$K$-theory give
the map of  complexes
\begin{equation}\label{CupProd}
p_1^*\eqref{K1Comp}\otimes p_2^*\eqref{K1Comp}\to \eqref{K2Comp} 
\end{equation}
over the cup product 
\begin{equation}\label{CupProd2}
p_1^*\sK_1\otimes p_2^*\sK_1\to \bar\sK_2.
\end{equation}

The augmentation $\epsilon_1:\sK_1\to \ker\beta$ is an isomorphism. 
The augmentation  $\epsilon_2:\bar\sK_2\to\ker\alpha$ is an injection, 
and the
cokernel is supported on $*\times *$, so
$\epsilon_2:\bar\sK_2\to\ker\alpha$ induces an isomorphism on $H^2$. 
Thus, the
complexes
\eqref{K1Comp} and \eqref{K2Comp} give rise to maps
\begin{gather*}
\delta_2:K_2(\C)\to H^2(X,\ker\alpha)=H^2(X,\bar\sK_2)=H^2(X,\sK_2)\\
\delta_1:\C^*=K_1(\C)\to H^1(E_0,\sK_1).
\end{gather*}
The compatibility of \eqref{CupProd} with \eqref{CupProd2} yields the 
commutativity of
the diagram
\[
\xymatrix{
\C^*\otimes\C^*\ar[r]^{\cup}\ar[d]_{\delta_1\otimes\delta_1}&K_2(\C)\ar[d]^
{\delta_2}\\
H^1(E_0,\sK_1)\otimes H^1(E_0,\sK_1)\ar[r]_-{p_1^*\cup p_2^*}& 
H^2(X,\sK_2).
}
\]
Since $L_v^\alg =\delta_1(v)$ for each $v\in\C^*$, we have
\[
p_1^*L_u^\alg\cup p_2^*L_{1-u}^\alg=\delta_2(\{u,1-u\})=0.
\]
\end{proof}

 The main
point of this section is that the Steinberg relation is {\em not} 
satisfied in
$\CH_0(E\times E')$ except in the case
$E=E'=E_0$. We first require the following lemma:

\begin{lem}\label{NonAlgLem} Let $s:\C\setminus\{0,1\}\to E\times E'$ be 
the analytic
map
$s(u)=(p(u),p'(1-u))$. Then
$s(\C\setminus\{0,1\})$ is not contained in any algebraic curve on 
$E\times E'$,
except in case
$E=E'=E_0$.
\end{lem}

\begin{proof} We first consider the case in which both $E$ and $E'$ are 
smooth
elliptic curves, $E=E_t$, $E'=E_{t'}$, where $t$ and $t'$ are in 
$\C^*$ and $|t|<1$, $|t'|<1$. We have the maps
\[
p:\C^*\to E,\ p':\C^*\to E',
\]
which are  group homomorphisms with $\ker p=t^\Z$, $\ker p'= 
t^{\prime\Z}$.

Suppose that $s(\C^*)$ is contained in an algebraic curve $D\subset 
E\times E'$.
For each $x\in E$, $(x\times E')\cap D$ is a finite set (possibly 
empty), hence,
for each $u\in\C\setminus\{0,1\}$, the set of points of $\C^*\times\C^*$ 
of the form
$(t^nu,1-t^nu)$ has finite image in $E\times E'$. Thus, for each $u$, 
there are
integers $n$, $m$ and $p$, depending on $u$, such that $n\neq m$ and
\begin{equation}\label{FinEq}
1-t^mu=t^{\prime p}(1-t^nu).
\end{equation}
Since there are uncountably many $u$, there is a single choice of $n$, 
$m$ and $p$
for which \eqref{FinEq} holds for uncountably many $u$. But then
\begin{equation}\label{FinEq2}
(t^{\prime p}t^n-t^m)u=1-t^{\prime p}.
\end{equation}
If $t^{\prime p}t^n-t^m=0$, then $|t'|=1$, contradicting the condition 
$|t'|<1$. If
$t^{\prime p}t^n-t^m\neq0$, then we can solve \eqref{FinEq2} for $u$, so
\eqref{FinEq} only holds for this single $u$, a contradiction.

If say $E'=E_0$, then $p':\C^*\to E'$ is injective, and we have the 
infinite set of
points $p'(1-t^nu)$ in the image of $s$, all lying over the single point 
$p(u)$.
\end{proof}

\begin{thm}\label{thm2.3} 
Let $E=E_t$, $E'=E_{t'}$, with at least one of $E$, $E'$ 
non-singular.
Then, for all $u$ outside a countable subset of $\C\setminus\{0,1\}$, 
$p(u)*p'(1-u)$ is not a torsion element in $F^2\CH_0(E\times E')$.
\end{thm}

\begin{proof} We first give the proof in case $E$ and $E'$ are both 
non-singular. For
a quasi-projective $\C$-scheme $X$, we let $S^nX$ denote the $n$th 
symmetric power of
$X$. For $X$ smooth, we have the map
\begin{align*}
\rho_n: S^nX(\C)\times S^nX(\C)&\to \CH_0(X)\\
(\sum_{i=1}^nx_i,\sum_{j=1}^ny_j)&\mapsto 
[\sum_{i=1}^nx_i-\sum_{j=1}^ny_j].
\end{align*}

For each integer $n\ge1$, we have the morphism
\begin{align*}
\phi_n:E\times E'&\to  S^{2n}(E\times E')\times S^{2n}(E\times E') \\
(x,y)&\mapsto(n(x,y)+n(0,0),n(x,0)+n(0,y)),
\end{align*}
By \cite[Theorem 1]{Roit},
$(\rho_{2n}\circ\phi_n)^{-1}(0)$ is a countable union of 
Zariski closed subsets  of $E\times E'$. 

On the other hand, since $p_g(E\times E')=1$, the Albanese kernel 
$F^2\CH_0(E\times
E')$ is ``infinite dimensional" \cite{Mumford}; in particular,  
$F^2\CH_0(E\times 
E')_\Q\neq0$.
Since $F^2\CH_0(E\times E')$ is generated by cycles of the form 
$p(u)*p(v)$, it
follows that $(\rho_{2n}\circ\phi_n)^{-1}(0)$ is a countable union of {\em
proper} closed subsets of $E\times E'$. If $D$ is a proper algebraic 
subset of
$E\times E'$, then, by Lemma~\ref{NonAlgLem}, $s^{-1}(D)$ is a proper 
analytic subset
of $\C\setminus\{0,1\}$, hence $s^{-1}(D)$ is countable. Thus, the set of 
$u\in
\C\setminus\{0,1\}$ such that $p(u)*p'(1-u)$ is torsion is countable, 
which completes
the proof in case both $E$ and $E'$ are non-singular.

If say $E'=E_0$, we use essentially the same proof. We let $X$ be the 
open subscheme
$E\times(E_0\setminus\{*\})$ of $E\times E_0$. We have the map
$\rho_n: S^nX(\C)\times S^nX(\C) \to \CH_0(E\times E_0)$ defined as above.
By \cite[Theorem 4.3]{LevWeib}, $(\rho_{2n}\circ\phi_n)^{-1}(0)$ is a 
countable union of  
closed subsets $D_i$ of $X$. By 
\cite{Srini}, we have the similar infinite dimensionality result for 
$\CH_0(E\times
E_0)$ as in  the smooth case, from which it follows that each $D_i$ is a 
proper
closed subset of 
$X$. Thus, the closure of each $D_i$ in $E\times E_0$ is a proper 
algebraic
subset of $E\times E_0$. The same argument as in the smooth case finishes 
the proof.
\end{proof}

\section{Indetectability} The zero-cycle $p(u)*p(1-u)$ is indetectable by 
cohomology
theories built on the sheaf $\sO_{E_\an\times E'_\an}^*$. We first 
consider the following
abstract situation.

Let $\Gamma_0(2)$ be the complex:
\begin{align*}
\Z[\C\setminus\{0,1\}]&\to \C^*\otimes\C^*\\
u&\mapsto u\otimes(1-u),
\end{align*}
with $\C^*\otimes\C^*$ in degree two.  

Let $X=E\times E'$, and let $\Gamma(2)_\an$ be a complex of sheaves on 
$X_\an$ with
the following properties:
\begin{equation}\label{Props}
\end{equation}
\begin{enumerate}
\item There is a group homomorphism $\cl:\CH_0(X)\to 
\H^4(X_\an,\Gamma(2)_\an)$.
\item There is a map in the derived category of sheaves
$D^b({\rm Sh}_{X_\an})$, $\rho:\sO^*_{X_\an}\otimes
\sO^*_{X_\an}[-2]\to
\Gamma(2)_\an$.
\item The composition 
\[
\C^*\otimes\C^*[-2]\to 
\sO^*_{X_\an}\otimes\sO^*_{X_\an}[-2]\to\Gamma(2)_\an
\]
extends to a map in $D^b({\rm Sh}_{X_\an})$, $\Gamma_0(2)\to 
\Gamma(2)_\an$.
\item The composition
\begin{multline*}
\Pic(X)\otimes \Pic(X)\cong H^1(X_\an,\sO^*_{X_\an})\otimes 
H^1(X_\an,\sO^*_{X_\an})\\
\xrightarrow{\cup}H^2(X_\an,\sO^*_{X_\an}\otimes \sO^*_{X_\an})
\xrightarrow{\rho}\H^4(X_\an,\Gamma(2)_\an)
\end{multline*}
agrees with the composition
\[
\Pic(X)\otimes
\Pic(X)\xrightarrow{\cup}\CH_0(X)\xrightarrow{\cl}\H^4(X_\an,\Gamma(2)_\an)
.
\]
\end{enumerate}
 
\begin{thm} \label{thm3.1}
Let $E=E_t$ and $E'=E_{t'}$, and let $\Gamma(2)_\an$ be a 
complex of
sheaves on $E_\an\times E'_\an$ satisfying the conditions \eqref{Props}. 
Then
$\cl(p(u)*p(1-u))=0$ for all
$u\in\C\setminus\{0,1\}$.
\end{thm}

\begin{proof} We give the proof in case both $E$ and $E'$ are 
non-singular; the
singular case is similar, but easier, and is left to the reader.

Since   
\[
p(u)*p(1-u)=[p_1^*c_1(L_u^\alg)]\cap[p_2^*c_1(L_{1-u}^\alg)],
\]
it follows from \eqref{Props}(4) that we need to show that
$\rho([L_u^\an]\cup[L_{1-u}^\an])=0$. The class $[L_u^\an]\in
H^1(E_\an,\sO^*_{E_\an})$ is the image of $[L_u]\in H^1(E_\an,\C^*)$ 
under the map of
sheaves $\C^*\to
\sO^*_{E_\an}$, and similarly for $L_{1-u}$ and $L_{1-u}^\an$. Thus, by
\eqref{Props}(3), it suffices to see that $p_1^*[L_u]\cup 
p_2^*[L_{1-u}]\in
H^2(E\times E',\C^*\otimes\C^*)$ vanishes in $\H^4(E\times 
E',\Gamma_0(2))$.

The $\Z$-covers $p:\C^*\to E=E_t$, $p':\C^*\to E'=E_{t'}$ give natural 
maps
\begin{gather}
\alpha:H^*(\Z,H^0(\C^*,\C^*))\to H^*(E_\an,\C^*), \notag\\
\beta:H^*(\Z,H^0(\C^*,\C^*))\to H^*(E'_\an,\C^*). \notag
\end{gather}
Similarly, the $\Z^2$-cover $p\times p':\C^*\times\C^*\to E\times E'$ 
gives the
natural map
\[
\gamma:\H^*(\Z^2,H^0(\C^*\times\C^*,\Gamma_0(2)))\to \H^*(E_\an\times 
E'_\an,\Gamma_0(2)).
\]
Letting $\iota:\C^*\otimes\C^*\to\Gamma_0(2)$ denote the natural 
inclusion, the maps
above are compatible with the respective cup products:
\[
\iota\circ(\alpha(a)\cup\beta(b))=\gamma\circ\iota(a\cup b).
\]

Each $v\in\C^*$ gives the corresponding homomorphism $v:\Z\to\C^*$, 
$v(n)=v^n$.
Since $[L_u]\in H^1(E_\an,\C^*)$ is $\alpha(u:\Z\to\C^*)$ and 
$[L_{1-u}]\in
H^1(E'_\an,\C^*)$ is $\beta(1-u:\Z\to\C^*)$, it suffices to show that
$\iota(p^*_1u\cup p_2^*(1-u))=0$ in 
$\H^4(\Z^2,\Gamma_0(2))$, where $p_1^*u,p_2^*(1-u):\Z^2\to\C^*$ are the 
respective homomorphisms
$(a,b)\mapsto u^a$, and $(a,b)\mapsto (1-u)^b$.

We have the spectral sequence
\[
E_2^{p,q}=H^p(\Z^2,H^q(\Gamma_0(2)))\Longrightarrow 
\H^{p+q}(\Z^2,\Gamma_0(2)).
\]
Since $\Z^2$ has cohomological dimension two, and since 
$H^q(\Gamma_0(2))=0$ for
$q\ne1,2$, it follows that the natural map
$\H^4(\Z^2,\Gamma_0(2))\to H^2(\Z^2,H^2(\Gamma_0(2)))$
is an isomorphism. Since $H^2(\Gamma_0(2))=K_2(\C)$, we need to show that 
the image
of $p_1^*u\cup p_2^*(1-u)$ in $H^2(\Z^2,K_2(\C))$ is zero.

By definition of the cup product in group cohomology, we have
\begin{align*}
[p_1^*u\cup p_2^*(1-u)]((a,b),(c,d))&= p_1^*u(a,b)\otimes 
p_2^*(1-u)(c-a,d-b)\\
&=
u^a\otimes (1-u)^{d-b},
\end{align*}
which clearly vanishes in $K_2(\C)$.

\end{proof}

\begin{ex}\label{ex2} In \cite{Bl}, S. Bloch defines a quotient
complex $\sB(2)$ of the analytic 
complex $ \sO^*_{X_\an}(1) \xrightarrow{2\pi i \otimes 1} \sO_{X_\an} 
\otimes
\sO^*_{X_\an}$ 
fulfilling 
$\sH^i(\sB(2)) = 0$ for $i \neq 1, 2$, 
$$\sH^1(\sB(2))= {\rm Im}\Big(r: K_{3, {\rm ind}}(\C) \to \C/\Z(2)\Big)
=:\Delta^*(1),$$ where $r$ is the regulator map,
and $\sH^2(\sB(2))= \sK_{2, {\rm an}}$. He shows
in the same article that $r(K_{3, {\rm ind}})(\C))= r(K_{3, {\rm
ind}})(\bar{\Q}))$, thus $\Delta^*(1)$ is a countable
subgroup of $\C/\Z(2)$, and also that $\sB(2)$ maps to the
complex $\Z(2) \to \sO_{X_\an} \to \Omega^1_{X_\an}$ which computes the
Deligne cohomology $H_{\sD}^*(X, 2)$ when $X$ is projective
smooth over $ \C$. In fact, the cycle map $\CH^2(X) \to
H^4_{\sD}(X, 2)$ is shown to factor through $H^4_{\sD}(X_{{\rm
an}}, \sB(2))$ (\cite{E}).
S. Bloch (\cite{B}) asked whether 
the cycle map $\CH^2(X) \to H^4_{\sD}(X_{{\rm
an}}, \sB(2))$  could possibly be injective. The
computations of this article show that it is not.
Indeed, by Lemma (1.3) of \cite{Bl}, the complex
$\Gamma_0(2)$ maps to the complex
$$\epsilon (\Z[\C\setminus \{0,1\}]) \to \C \otimes \C^*,$$
where $\epsilon$ is defined via the dilogarithm function
$$\epsilon(a)= [\log (1-a) \otimes a] + [2\pi i \otimes {\rm
exp} \Big( \frac{-1}{2\pi i} \int_0^a \log (1-t) \frac{dt}{t}
\Big)],$$ 
and the latter complex maps to 
$$\sB(2)_{X}: 
\sO^*_{X_\an}(1) \xrightarrow{2\pi i \otimes 1} \sO_{X_\an} \otimes
\sO^*_{X_\an}/(\epsilon  (\Z[\C\setminus \{0,1\}])$$
for $X= \Spec \C$.
Let us take $\Gamma(2)_\an= \sB(2)$. We now
verify the conditions \ref{Props}.
The condition 1 is given by \cite{E}. Indeed, one computes
the Leray spectral sequence associated to
$\alpha: X_\an \to X_{{\rm Zar}}$ and the first term entering
$H^4(\sB(2))$ is $$E^{2,2}= H^2_{\rm Zar}(R\alpha_* \sB(2))=
H^2(\sK_{2, \Z}),$$ where $\sK_{2, \Z} := {\rm Ker}\Big(
\alpha_*\sK_{2,\an} \xrightarrow{d\log \wedge d\log} H^2(\C/\Z(2))\Big).$
Then the cycle map 
$\cl$ is induced by $\sK_2 \to \sK_{2, \Z}$ on $X_{{\rm Zar}}$,
which is obviously compatible with the product in ${\rm Pic}$.
Thus we have 4. 
We have already discussed 2 and 3. Hence we can apply 
Theorem \ref{thm2.3} to take a 0-cycle $p(u)*p(1-u)$
on $E \times E'$ where both $E$ and $E'$ are smooth elliptic
curves which does not die in the Chow group $\CH_0(E \times
E')$, whereas it dies by Theorem \ref{thm3.1} in $\H^4(\sB(2))$.

In \cite{Li}, S. Lichtenbaum constructs an \'etale version
$\Gamma(2)$ of S. Bloch's analytic complex $\sB(2)$, the
cohomology of which contains $\CH^2(X)$. This contrasts
with the examples discussed above. 

Over a $p-$adic field, 
W. Raskind and M. Spie{\ss} (\cite{RS}) show that
the Albanese
kernel modulo $n$ of a product of two Tate elliptic curves
is dominated by $K_2(k)/n$. This result is not immediately
comparable to ours, but is obviously related.

\end{ex}
\section{The Relative Situation}

In this section, we study the cycles constructed in section
2 on $X=E \times E_0$, where as there,
$E$ is smooth, and $E_0$ is a nodal curve.
Let $\nu=1\times q : E \times \P^1 \to X$ be the normalization.
We define
\begin{gather}
\bar{\sK_2}= {\rm Ker} \Big(\nu_*\sK_2 
\xrightarrow{|_{E \times 0} - 
|_{E \times \infty}} \sK_2|_{E }\Big)
\end{gather}

\begin{lem} \label{lem4.1} One has
\[
\CH^2(X)= H^2(X, \bar{\sK_2}),
\]
and the Chow group $\CH_0(X)$
fits into an exact sequence
\[
0 \to H^1(E, \sK_2) \xrightarrow{\gamma} \CH_0(X)
\xrightarrow{\nu^*} \CH_0(E \times \P^1)={\rm Pic}(E)\otimes
{\rm Pic}(\P^1) \to 0.
\]
Moreover, the  map $\gamma$ is defined by
\[
\gamma(\sum_{x \in E^{(1)}} x \otimes \lambda_x) = \sum_{x \in
E^{(1)}} (x, p_0(\lambda_x))-(x,0).
\]
\end{lem}
\begin{proof} 
The map $\nu^*: \sK_2 \to \bar{\sK_2}$ is obviously surjective, and by the
Gersten resolution on the smooth points of $X$, 
the kernel is supported in codimension 1.
Thus $\nu^*$ induces an isomorphism on $H^2$. 

On the other hand, 
$$H^1(E \times \P^1, \sK_2)= H^1(E, \sK_2)
\oplus H^0(E, \sK_1) \cup c_1(\sO(1)).$$ The term $H^1(E, \sK_2)$
maps to $0 \in H^1(E, \sK_2)$ via the difference of the
restrictions to $E \times 0$ and $E \times \infty$, while
$c_1(\sO(1))$ restricts to 0 to either $E\times 0$ or $E \times
\infty$. This shows the long exact
sequence associated to the short one defining $\bar{\sK_2}$.

Finally, the value $\gamma(x\otimes \lambda_x)$ of the 
map is given by the boundary morphism
$\C^* \to H^1(X, \sO^*_X)$ induced by the normalization sequence
$$0 \to \sO^*_X \to  q_*\sO^*_{\P^1} \xrightarrow{|_0 -
|_{\infty}} \C^* \to 0$$
on the right argument $\lambda_x$. The formula for $\gamma$ thus follows 
from
Lemma~\ref{C1Comp}.  
\end{proof}

Let ${\rm Nm}: H^1(E, \sK_2) \to \C^*$ be the norm map defined by
\begin{gather}\label{Nm}
{\rm Nm}\Big(\sum_{x \in E^{(1)}}x \otimes \lambda_x \Big) =
\prod_{x \in E^{(1)}} \lambda_x.
\end{gather}
We set
\begin{gather}\label{dfnV}
V(E) ={\rm Ker} {\rm Nm}.
\end{gather}
One has
\begin{lem} \label{lem4.2}
$F^2\CH_0(X) = \gamma\Big(V(E)\Big).$
\end{lem}
\begin{proof}
By the definition given in \S\ref{AlbKer}, $F^2\CH_0(X)$ is generated by 
the
expressions 
$[(x,y)]-[(x,0)]-[(0,y)]+[(0,0)]$, with $x\in E(\C)$ and $y\in 
E_0(\C)\setminus\{*\}$.
By the formula for
$\gamma$ given in Lemma~\ref{lem4.1}, this expression is $\gamma(x\otimes 
y-0\otimes
y)$, after identifying $y\in\C^*$ with $p_0(y)\in E_0(\C)$. Clearly 
$V(E)$ is
generated by  the elements of $H^1(E,\sK_2)$ of the form $x\otimes 
y-0\otimes y$,
whence the lemma.
\end{proof}

Next we want to map $\CH_0(X)$ to a relative version of S.
Bloch's analytic motivic cohomology. So we define
\begin{gather}
\bar{\sB}(2) = {\rm Ker} \Big(\nu_*\sB(2) 
\xrightarrow{|_{E \times 0} - 
|_{E \times \infty}} \sB(2)|_{E }\Big)
\end{gather}
In particular, $\bar{\sB}(2)$ is an extension of 
$$\bar{\sK}_{2, \an}= {\rm Ker} \Big(\nu_*\sK_{2,\an} 
\xrightarrow{|_{E \times 0} - 
|_{E \times \infty}} \sK_{2, \an}|_{E }\Big)$$
placed in degree 2, by $\Delta^*(1)$, placed in degree 1.
In other words, $\sB(2)$ is the pull-back of $\bar{\sB}(2)$  via
the map $ \nu^*:\sK_{2, \an} \to \bar{\sK}_{2, \an}$, and in
particular, it receives the complex $\Gamma_0(2)$ as explained
in the example \ref{ex2}.

Considering again the Leray spectral sequence attached
to the identity $\alpha: X_\an \to X_{{\rm zar}}$, we 
see that
\begin{gather}
\bar{\sK}_{2, \Z}  := {\rm Ker}\Big(\alpha_* \bar{\sK}_{2, \an} 
\to \sH^2(\C/\Z(2))\Big)
\end{gather}
receives $\bar{\sK_2}$ and that the first map of the spectral
sequence is then
\begin{gather}
H^2(X, \bar{\sK}_{2, \Z}) \to \H^4(X_\an, \bar{\sB}(2)).
\end{gather}
In conclusion, we have shown
\begin{lem}
One has a cycle map
$$\psi_X: \CH_0(X) \to \H^4(X_\an, \bar{\sB}(2))$$
compatible with the cycle map
$$\psi_{E \times \P^1}:
\CH_0(E \times \P^1) \to \H^4((E \times \P^1)_\an, \sB(2))$$
on the normalization. 
Moreover, $\psi_X$ fulfills the conditions described in \ref{Props}.
\end{lem}
\begin{proof}
We just have to verify the condition 4 of \ref{Props}. 
{From} the normalization sequence
$$0 \to \sO^*_X \to \nu_*\sO^*_{E\times \P^1}  
\xrightarrow{|_{E \times 0} - 
|_{E \times \infty}} \sO^*_E \to 0,$$
one has a natural map
$$\sO^*_{X_\an} \otimes \sO^*_{X_\an} \to \bar{\sK}_{2,\an}$$
which obviously fulfills \ref{Props} 4. 
\end{proof}

Now we can apply Theorem \ref{thm3.1} to conclude
\begin{thm}\label{SingVanThm}
The 0-cycles defined by the Steinberg curve on $E\times E_0$ die
in the analytic motivic cohomology $\H^4(X_\an, \bar{\sB}(2))$.
\end{thm}

Let $K$ be a subfield of $\C$.
We next consider for any algebraic variety $Z$ defined over $K$,
the cycle map with values in the absolute 
Hodge cohomology
\begin{gather}
H^m(Z, \sK_2) \xrightarrow{d\log \wedge d\log}H^m(Z,
\Omega^2_{Z/\Q}) 
\end{gather}
induced by the absolute $d\log$ map
\begin{gather}
\sO^*_Z \xrightarrow{d\log} \Omega^1_{Z/\Q}.
\end{gather}
This cycle map is obviously compatible with the map
$\gamma$, and with extension of scalars. 

Let $E\to\Spec K$ be an elliptic curve over a subfield $K$ of $\C$. We 
have the exact
sheaf sequence 
\[
0\to\sO_E\otimes\Omega^1_{K/\Q}\to\Omega^1_{E/\Q}\to\Omega^1_{E/K}\to0,
\]
which induces a two-term filtration $F^*\Omega^2_{E/\Q}$ of 
$\Omega^2_{E/\Q}$ with
$F^2\Omega^2_{E/\Q}=\sO_E\otimes\Omega^2_{K/\Q}$. This gives us the 
natural maps
\begin{gather*}
\gamma_1:H^*(E,\sO_E)\otimes\Omega^1_{K/\Q}\to H^*(E,\Omega^1_{E/\Q})\\
\gamma_2:H^*(E,\sO_E)\otimes\Omega^2_{K/\Q}\to H^*(E,\Omega^2_{E/\Q}).
\end{gather*}

We have the norm map $\Nm:H^1(E,\sK_2)\to H^0(K,\sK_1)=K^*$ as
in \ref{Nm}, but over $K$; we let 
$V(E)\subset
H^1(E,\sK_2)$ be the kernel of $\Nm$ (see \eqref{dfnV}).

\begin{lem}\label{AHCompLem} Let $K$ be an algebraically closed subfield 
of $\C$,
$E\to\Spec K$ an elliptic curve over $K$.  Then the cycle map with values 
in
absolute Hodge cohomology maps $V(X)$ to the subgroup $\gamma_2[H^1(E, 
\sO_E)
\otimes\Omega^2_{E/\Q}]$ of $H^1(E,\Omega^2_{E/\Q})$.
\end{lem}

\begin{proof} The kernel of the composition
\[
\Pic(E)=H^1(E,\sK_1)\xrightarrow{d\log}H^1(E,\Omega^1_{E/\Q})\to
H^1(E,\Omega^1_{E/K})\cong K
\]
is the composition
\[
\Pic(E)\xrightarrow{\deg}\Z\subset K,
\]
hence the $d\log$ map sends $\Pic^0(E)$ to the subgroup
$\gamma_1[H^1(E,\sO_E)\otimes\Omega^1_{K/\Q}]$ of 
$H^1(E,\Omega^1_{E/\Q})$.

Take $\tau\in\Pic^0(E)$, $u\in H^0(E,\sK_1)=K^*$, and let $\xi=\tau\cup 
u\in
H^1(E,\sK_2)$. Then
\[
d\log(\xi)=d\log(\tau)\cup d\log(u).
\]
Since $d\log:K^*\to\Omega^1_{K/\Q}$ is just the absolute $d\log$ map, we 
see that
$d\log(\xi)$ lands in the image of the cup product map
\[
[H^1(E,\sO_E)\otimes\Omega^1_{K/\Q}]\otimes \Omega^1_{K/\Q}\to 
H^1(E,\Omega^2_{E/\Q}),
\]
which is $\gamma_2(H^1(E,\sO_E)\otimes\Omega^2_{K/\Q})$.

Since $K$ is algebraically closed, the cup product $\Pic(E)\otimes K^*\to
H^1(E,\sK_2)$ is surjective, from which one sees that the cup product maps
$\Pic^0(E)\otimes K^*$ onto $V(E)$. Combining this with the computation 
above
completes the proof.
\end{proof}

{From} the surjectivity of the cup product $\Pic^0(E)\otimes K^*\to V(E)$ 
for $K$
algebraically closed, we see that the injection $H^1(E,\sK_2)\to 
\CH_0(X)$ sends
$V(E)$ isomorphically onto $F^2\CH_0(X)$.

Let $K$ be a subfield of $\C$. We say that an element $\xi$ of $\CH_0(X)$ 
is {\em
defined over $K$} if there is an $K$-scheme $X^0$, an element $\xi^0$ of
$\CH_0(X^0)$ and an isomorphism $\alpha:X^0_\C\to X$ such that
$\xi=\alpha_*(\xi^0_\C)$. From Lemma~\ref{AHCompLem} and the 
compatibility of $d\log$
with extension of scalars, we have

\begin{lem}\label{AHVanLem} Take $K=\C$, and let $\xi$ be an element of
$F^2\CH_0(X)=V(E)$. If
$\xi$ is defined over a field of transcendence degree one over $\Q$, then 
$\xi$
vanishes under the cycle map to absolute Hodge cohomology.
\end{lem}

\begin{cor} 
If $E$ is an elliptic curve with complex multiplication, then
there are non-torsion cycles $\xi \in F^2\CH_0(X)$ dying in the
analytic motivic cohomology as well as in absolute Hodge
cohomology. \end{cor}

\begin{proof} By the remark above, we may replace $F^2\CH_0(X)$ with 
$V(E)$. Let
$\bar E$ be a model for $E$, with equation
$y^2=4x^3-ax-b$ defined over a number field
$K\subset\C$. Let $\omega=\frac{dx}{y}$ be the standard global one-form 
on $\bar E$.

Choosing an isomorphism $\bar E_\C\cong E_\C$ defines the period lattice
$L_\omega\subset\C$ for $\omega$. Choose a basis for
$L_\omega$ of the form $\{\Omega,\tau\Omega\}$, and let $t=e^{2\pi 
i\tau}$. Let
\[
\sP:\C\to \C\P^1
\]
be the Weierstra\ss\ $P$-function for the lattice $L_\omega$.

The map $\times\Omega^{-1}:\C\to\C$ gives rise to the
isomorphism of Riemann surfaces
$\alpha_\an:\bar E_\C^\an\to E_t^\an$ making the diagram 
\[
\xymatrix{
\C\ar[r]^{\times\Omega^{-1}}\ar[dd]_{(\sP,\sP')}&\C\ar[d]^{\text{exp}}\\
&\C^*\ar[d]^p\\
\bar E_\C^\an\ar[r]_{\alpha_\an}&E_t^\an}
\]
commute, i.e.,
\[
p(u)=\alpha_\an(\sP(\frac{\Omega}{2\pi i}\text{log}u), 
\sP'(\frac{\Omega}{2\pi
i}\text{log}u)).
\]
We let
\[
\alpha:\bar E_\C\to E_t
\]
be the corresponding isomorphism of algebraic elliptic curves over $\C$.

By \cite{Be}, th\'eor\`eme 1,  $\sP(\frac{\Omega}{2\pi i}\text{log}u)$
has transcendence degree 1 over $\bar\Q$ for all 
$u\in \N$, $u\ge2$. (We thank Y. Andr\'e for 
giving us this reference). Fix a $u\ge2$,  let $K$ be the algebraic 
closure of the
field $\Q(\sP(\frac{\Omega}{2\pi i}\text{log}u))$, and let $x\in\bar 
E(K)$ be the
point $(\sP(\frac{\Omega}{2\pi i}\text{log}u), \sP'(\frac{\Omega}{2\pi
i}\text{log}u))$. Then $x$ is a generic point of $\bar E$ over $\bar\Q$.

We take 
\[
\xi:= p(u)*p(1-u).
\]
By construction, $\xi=\alpha(\xi_K\times_K\C)$, where $\xi_K\in H^1(\bar 
E,\sK_2)$ is
the element $[(x)-(0)]\cup[1-u]$. Here  $[(x)-(0)]$ denotes the class in
$\Pic(E)=H^1(\bar E,\sK_1)$, and $[1-u]$ denotes the class in $H^0(\bar
E,\sK_1)=K^*$. Since $K$ has transcendence degree one over $\bar\Q$,  the 
class
of
$\xi$ in the absolute Hodge cohomology of $E$ vanishes, by 
Lemma~\ref{AHVanLem}. By
Theorem~\ref{SingVanThm}, $\xi$ dies  in the analytic  motivic cohomology 
of $E$ as
well. It remains to show that $\xi$ is a non-torsion element of 
$H^1(E_K,\sK_2)$.

We give an analytic proof of this using the regulator map with
values in Deligne-Beilinson cohomology.

Let $Y$ be a smooth projective surface over $\C$, and let $\NS(Y)$ denote 
the
N\'eron-Severi group of divisors modulo homological equivalence.
Then Hodge theory implies that
$$ \NS(Y)= \{(z, \varphi) \in (H^2(Y_\an, \Z(1)) \times
F^1H^2(Y_\an, \C)), z\otimes \C= \varphi\},$$
and that
$$\NS(Y) \cap F^2H^2_{DR}(Y) =\emptyset.$$
We note that the map $\Pic(Y)\otimes\C^*\to H^3_\sD(Y,\Z(2))$ induced by 
the cup
product in Deligne cohomology factors through $\NS(Y)\otimes\C^*$, and 
that the
induced map $\iota:\NS(Y)\otimes\C^*\to H^3_\sD(Y,\Z(2))$ is injective. 
Indeed, 
$$
H^3_\sD(Y,\Z(2))= H^2(Y_\an, \C/\Z(2))/F^2.$$

Now take $Y=E\times E$, and
let $U\subset E$ be the complement of a non-empty finite set $\Sigma$ of 
points of
$E$. Let $[E\times0]$ be the class of $E\times 0$ in $\NS(Y)$, and let 
$\gamma:\C^*\to
\NS(Y)\otimes\C^*$ be the map $\gamma(v)=[E\times0]\otimes v$. Let
\[
\iota_U:\NS(Y)\otimes\C^*\to H^3_\sD(E\times U,\Z(2))
\]
be the composition of $\iota$ with the restriction map 
$H^3_\sD(Y,\Z(2))\to
H^3_\sD(E\times U,\Z(2))$.  We claim that the sequence
\[
\C^*\xrightarrow{\gamma}\NS(Y)\otimes\C^*\xrightarrow{\iota_U}
H^3_\sD(E\times U,\Z(2))
 \]
is exact. Indeed, we have the localization sequence
\[
\oplus_{s\in\Sigma}H^1_\sD(E\times s,\Z(1))\xrightarrow{\oplus_s\iota_s}
H^3_\sD(Y,\Z(2))\to H^3_\sD(E\times U,\Z(2))\to,
\]
the isomorphism $H^1_\sD(E\times s,\Z(1))\cong\C^*$ and the identity
\[
\iota_s(v)=\gamma(v),\ v\in\C^*,
\]
which proves our claim.

In particular, let $[\Xi]=[\Delta-\{0\}\times E]\otimes v$, where 
$\Delta$ is 
the
diagonal,  $v$ is an element of $\C^*$ which is not a root of unity, and
$[\Delta-\{0\}\times E]$ is the class in $\NS(Y)$. Since
$[\Delta-\{0\}\times E]$ is not torsion in $\NS(Y)/[E\times
\{0\}]$, we see that 
$[\Xi]$ has
non-torsion image $[\Xi_{\C(E)}]$ in
\[
H^3_\sD(E\times_\C\C(E),\Z(2)):=\lim_{\substack{\to\\\0\neq U\subset
E}}H^3_\sD(E\times U,\Z(2)),
\]
where the limit is over non-empty Zariski open subsets $U$ of $E$.

Let $\Xi$ be the image of $(\Delta-0\times E)\otimes v$ in 
$H^1(Y,\sK_2)$. Then
$[\Xi]$ is the image of $\Xi$ under the regulator map $H^1(Y,\sK_2)\to 
H^3_\sD(Y,\Z(2))$. Similarly, letting $\Xi_{\C(E)}$ be the pull-back of 
$\Xi$ to
$E\times_\C\C(E)$, $[\Xi_{\C(E)}]$ is the image of $\Xi_{\C(E)}$ under 
the 
regulator map $H^1(E\times_\C\C(E),\sK_2)\to 
H^3_\sD(E\times_\C\C(E),\Z(2))$. Thus,  $\Xi_{\C(E)}$ is a non-torsion 
element of 
$H^1(E\times_\C\C(E),\sK_2)$ for each non-torsion element $v\in\C^*$.

Let $\bar\Delta$ be the diagonal in $\bar E\times\bar E$, let $\bar\xi$ 
be the
image of  $(\bar\Delta-0\times \bar E)\otimes(1-u)$ in $H^1(E,\sK_2)$, 
and let
$\bar\xi_{\bar\Q(E)}$ be the image of $\bar\xi$ in
$H^1(\bar E\times_{\bar\Q}\bar\Q(\bar E),\sK_2)$. Clearly, after
choosing a complex embedding $\bar{\Q} \subset \C$,  $\Xi_{\C(E)}$ 
(for $v=1-u$)
is the image of $\bar\xi_{\bar\Q(E)}$ under the extension of scalars 
$\bar\Q(\bar
E)\to \C(\bar E)\cong\C(E)$, hence $\bar\xi_{\bar\Q(E)}$ is a non-torsion 
element of
$H^1(\bar E\times_{\bar\Q}\bar\Q(\bar E),\sK_2)$.

Since $x$ is a geometric generic point of $\bar{E}$ over
$\bar{\Q}$, there is an embedding   $\sigma:\bar\Q( E)\to\C$ such that
$x:\Spec\C\to\bar E$ is the composition $\Spec\C\to\Spec\bar\Q(E)\to \bar 
E$. Thus,
$\xi$ is the image of $\bar\xi$ under $(\id\times x)^*:H^1(\bar 
E\times_{\bar{\Q}}\bar
E,\sK_2)\to H^1(E,\sK_2)$, and hence $\xi$ is the image of 
$\bar\xi_{\bar\Q(E)}$
under the map $\id\times\sigma_*:H^1(\bar E\times_{\bar\Q}\bar\Q(\bar 
E),\sK_2)\to
H^1(E,\sK_2)$ induced by the extension of scalars $\sigma$.

Since the kernel of $\id\times \sigma_*$ is torsion, it
follows that $\xi$ is a non-torsion element of $H^1(E,\sK_2)$, as desired.
\end{proof}

\begin{rem}
Going back to $X=E \times E'$, where both elliptic curves are
smooth, we are lacking the transcendence theorem  which would 
force the existence of a cycle $0\neq \xi= p(u) * p(1-u)\in
F^2\CH_0(X)$  dying both in $\H^4(X, \sB(2))$ and in
absolute Hodge cohomology.
\end{rem}


\begin{thebibliography}{99}
\bibitem{Be} Bertrand, D.: Valeurs de fonctions th\^{e}ta
et hauteurs $p$-adiques, Progress in Math. {\bf 22} (1982), 1-11, 
Birkh\"auser Verlag.
\bibitem{SriniBis} Biswas, J.; Srinivas V.: Chow ring of a
singular surface, appendix to ``Roitmann theorem for singular
projective varieties'', preprint 1995, to appear in Compositio.
\bibitem{Bl} Bloch, S.: Applications of the dilogarithm function
in algebraic $K$-theory and algebraic geometry, Int. Symp. on
Alg. Geom., Kyoto (1977), 2036-2060.
\bibitem{B} Bloch, S.: letter to H. Esnault, March 30, 1988.
\bibitem{E} Esnault, H.: A note on the cycle map, J. reine
angew. Math. {\bf 411} (1990), 51-65.
\bibitem{LevWeib} Levine, M.; Weibel, C.: Zero-cycles and complete 
intersections on 
affine surfaces, J. Reine u. Ang. Math. {\bf 359}(1985) 106-120.
\bibitem{Lev} Levine, M.:  Bloch's formula for a singular surface, 
Topology {\bf
124} No. 2(1985) 165-174.
\bibitem{Li} Lichtenbaum, S.: The construction of weight-two
motivic cohomology, Invent. math. {\bf 88} (1987), 183-215.

\bibitem{Mumford}
Mumford, D.: Rational equivalence of $0$-cycles on surfaces, J. Math.
Kyoto Univ. {\bf 9} (1968) 195-204.

\bibitem{M}
Mumford, D.: Abelian Varieties, Oxford University Press (1970).

\bibitem{RS} Raskind, W.; Spie{\ss}, M.: Milnor $K$-groups
and zero-cycles on products of curves over ${p}$-adic fields,
preprint 31 pages.

\bibitem{Roit} Roitman, A. A.:
$\Gamma$-equivalence of zero-dimensional cycles.  
Mat. Sb. (N.S.) {\bf 86}(128) (1971) 557-570. 

\bibitem{Se} Serre, J.-P.: G\'eom\'etrie alg\'ebrique et
g\'eom\'etrie analytique, Ann. Inst. Fourier {\bf 6} (1956),
1-42. 

\bibitem{Srini} Srinivas, V.: Zero cycles on singular avrieties,
preprint 1998, 26 pages, appears in the Proceedings of the Banff
conference on algebraic cycles.
 
\end{thebibliography}
\end{document}